\documentclass[12pt]{article}
\usepackage{latexsym,amsfonts,amssymb}
\usepackage[dvips]{color}
\usepackage{colordvi,multicol}

\def \cal{\mathcal}

\newtheorem{thm}{Theorem}[section]
\newtheorem{cor}[thm]{Corollary}
\newtheorem{lem}[thm]{Lemma}
\newtheorem{pro}[thm]{Proposition}
\newtheorem{defi}[thm]{Definition}
\newtheorem{rem}[thm]{Remark}
\newtheorem{exa}[thm]{Example}

\begin{document}

%

\centerline{\Large\bf New results on Hunt's hypothesis (H)} \vskip
0.2cm \centerline{\Large\bf  for L\'{e}vy processes}

\vskip 1.6cm \centerline{Ze-Chun Hu} \centerline{\small Department
of Mathematics, Nanjing University, Nanjing, 210093, China}
\centerline{\small E-mail: huzc@nju.edu.cn}
\vskip 0.7cm \centerline{Wei Sun} \centerline{\small Department of
Mathematics and Statistics, Concordia University,}
\centerline{\small Montreal, H3G 1M8, Canada} \centerline{\small
E-mail: wei.sun@concordia.ca}
\vskip 0.7cm \centerline{Jing Zhang} \centerline{\small Department of
Mathematics and Statistics, Concordia University,}
\centerline{\small Montreal, H3G 1M8, Canada} \centerline{\small
E-mail: waangel520@gmail.com} \vskip 1.4cm

\vskip 0.5cm \noindent{\bf Abstract}\quad In this paper, we
present new results on Hunt's hypothesis (H) for L\'{e}vy
processes. We start with a comparison result on L\'{e}vy processes
which implies that big jumps have no effect on the validity of
(H). Based on this result and the Kanda-Forst-Rao theorem, we give
examples of subordinators satisfying (H). Afterwards we give a new
necessary and sufficient condition for (H) and obtain an extended
Kanda-Forst-Rao theorem. By virtue of this theorem, we give a new
class of L\'{e}vy processes satisfying (H). Finally, we construct
a type of subordinators that does not satisfy Rao's condition.

\smallskip

\noindent {\bf Keywords}\quad Hunt's hypothesis (H), Getoor's conjecture,  L\'{e}vy process, subordinator.

\smallskip

\noindent {\bf Mathematics Subject Classification (2010)}\quad
Primary: 60J45; Secondary: 60G51

\tableofcontents

\section{Introduction}
Let $X$ be a time-homogeneous Markov process. Hunt's hypothesis
(H) says that ``every semipolar set of $X$ is polar". This
hypothesis plays a crucial role in the potential theory of (dual)
Markov processes. To illustrate its importance, let us recall some
potential-theoretic principles.

Suppose that $E$ is a locally compact space with a countable base.
Let $(X, P^x)$ and $(\hat{X},\hat P^x)$ be a pair of dual standard
Markov processes on $E$ as described in Blumenthal and Getoor \cite[VI]{BG68}. Denote by
${\cal B}^n$ the family of all nearly Borel measurable subsets of
$E$. For $D\subset E$, we define the first hitting time of $D$ by
$$
\sigma_D:=\inf\{t>0:X_t\in D\}.
$$
A set $D\subset E$ is called polar (respectively, essentially
polar) if there exists a set $C\in {\cal B}^n$ such that $D\subset
C$ and $P^x(\sigma_C<\infty)=0$ for every $x\in E$ (respectively,
almost every $x\in E$ with respect to the reference measure). $D$
is called a thin set if there exists a set $C\in {\cal B}^n$ such
that $D\subset C$ and {$P^x(\sigma_C=0)=0$} for every $x\in E$.
$D$ is called semipolar if $D\subset\bigcup_{n=1}^{\infty}D_n$ for
some thin sets $\{D_n\}_{n=1}^{\infty}$.

 Denote by $E^x$ the expectation with respect to $P^x$.  Let $\alpha>0$. A finite $\alpha$-excessive function $f$ on $E$ is called a regular potential provided that
$E^x\{e^{-\alpha T_n}f(X_{T_n})\}\rightarrow E^x\{e^{-\alpha
T}f(X_{T})\}$  for $x\in E$ whenever $\{T_n\}$ is an increasing
sequence of stopping times with limit $T$. Denote by
$(U^{\alpha})_{\alpha>0}$ the resolvent operators for $X$.

\begin{itemize}
\item {\bf Bounded positivity principle $(P^*_{\alpha})$}: If $\nu$ is a finite signed measure
such that  $U^{\alpha}\nu$ is bounded, then $\nu U^{\alpha}\nu\geq 0$, where $\nu U^{\alpha}\nu:=\int_EU^{\alpha}\nu(x)\nu(dx)$.

\item {\bf Bounded energy principle $(E^*_{\alpha})$}: If $\nu$ is a finite measure with compact support
such that  $U^{\alpha}\nu$ is bounded, then $\nu$ does not charge
semipolar sets.

\item {\bf Bounded maximum principle $(M^*_{\alpha})$}: If $\nu$ is a finite measure
with compact support $K$ such that $U^{\alpha}\nu$ is bounded, then
$\sup\{U^{\alpha}\nu(x):x\in E\}=\sup\{U^{\alpha}\nu(x):x\in K\}$.

\item {\bf Bounded regularity principle $(R^*_{\alpha})$}: If $\nu$ is a finite measure
with compact support such that $U^{\alpha}\nu$ is bounded, then
$U^{\alpha}\nu$ is regular.

\item {\bf Polarity principle} ({\bf Hunt's hypothesis (H)}): Every semipolar set is polar.

\end{itemize}

\begin{pro}\label{thm1}
Assume that all 1-excessive (equivalently, all $\alpha$-excessive, $\alpha>0$) functions are lower
semicontinuous. Then
$$(P^*_{\alpha})\Leftrightarrow (E^*_{\alpha})
\Leftrightarrow(M^*_{\alpha})\Leftrightarrow(R^*_{\alpha})\Leftrightarrow ({\rm H}).
$$
\end{pro}
{\bf Proof.} $(R^*_{\alpha})\Leftrightarrow ({\rm H})$ is proved in Blumenthal and Getoor \cite{BG68} and $(M^*_{\alpha})\Leftrightarrow ({\rm H})$ is proved in Blumenthal and Getoor \cite{BG70}.
$(P^*_{\alpha})\Rightarrow (M^*_{\alpha})$ is proved in Rao \cite{R77} and
$(M^*_{\alpha})\Rightarrow (P^*_{\alpha})$ is proved in Fitzsimmons \cite{Fi90}.
By \cite[Propsition (2.1)]{BG70}, $(E^*_{\alpha})
\Rightarrow (M^*_{\alpha})$. By \cite[Proposition (5.1)]{BG70} and the equivalence of $(M^*_{\alpha})$ and $({\rm H})$, $(M^*_{\alpha})
\Rightarrow (E^*_{\alpha})$.\hfill\fbox

Hunt's hypothesis (H) is also equivalent to some other important
properties of Markov processes. For example, Blumenthal and Getoor
\cite[Proposition (4.1)]{BG70} and Glover \cite[Theorem
(2.2)]{G83} showed that (H) holds if and only if the fine and
cofine topologies differ by polar sets; Fitzsimmons and Kanda
\cite{FK} showed that (H) is equivalent to the dichotomy of
capacity.

In spite of its importance, (H) has been verified only in some
special situations.  Some forty years ago, Getoor conjectured that
essentially all L\'{e}vy processes satisfy (H).

From now on  we let $(\Omega,{\cal F},P)$ be a probability space
and $X=(X_t)_{t\ge 0}$ be an $\mathbf{R}^n$-valued L\'{e}vy
process on $(\Omega,{\cal F},P)$ with L\'{e}vy-Khintchine exponent
$\psi$, i.e.,
\begin{eqnarray*}
E[\exp\{i\langle z,X_t\rangle\}]=\exp\{-t\psi(z)\},\  z\in
\mathbf{R}^n,t\ge 0,
\end{eqnarray*}
where $E$ denotes the expectation {with respect to} $P$ and $\langle\cdot,\cdot\rangle$ denotes the Euclidean inner product of $\mathbf{R}^n$. The classical L\'{e}vy-Khintchine formula tells us that
\begin{eqnarray*}
\psi(z)=i\langle a,z\rangle+\frac{1}{2}\langle
z,Qz\rangle+\int_{\mathbf{R}^n} \left(1-e^{i\langle
z,x\rangle}+i\langle z,x\rangle 1_{\{|x|<1\}}\right)\mu(dx),
\end{eqnarray*}
where $a\in \mathbf{R}^n,Q$ is a symmetric nonnegative definite
$n\times n$ matrix, and $\mu$ is a measure (called the L\'evy
measure) on $\mathbf{R}^n\backslash\{0\}$ satisfying
$\int_{\mathbf{R}^n\backslash\{0\}} (1\wedge
|x|^2)\mu(dx)<\infty$. Hereafter, we use Re$(\psi)$ and Im$(\psi)$
to denote the real and imaginary parts of $\psi$, respectively,
and use $(a,Q,\mu)$ to denote $\psi$.

Let us recall some important results obtained so far for Getoor's
conjecture. When $n=1$, Kesten \cite{Ke69} {(cf. also Bretagnolle
\cite{Br71}) showed   that if $X$ is not a compound Poisson
process, then  every $\{x\}$ is non-polar} if and only if
\begin{eqnarray*}\label{Ke69-a}
\int_0^{\infty}\mbox{Re}([1+\psi(z)]^{-1})dz<\infty.
\end{eqnarray*}
Port and Stone \cite{PS69} proved that for the asymmetric
Cauchy process on the line every $x$ is regular for $\{x\}$, and thus (H) holds in this case. Further, Blumenthal and Getoor
\cite{BG70} showed that all  stable processes with index
$\alpha\in (0,2)$ on the line satisfy (H).

 Kanda \cite{Ka76} and Forst \cite{F75} proved that (H)
holds if $X$ has bounded continuous transition densities (with respect to the Lebesgue measure $dx$) and
the L\'{e}vy-Khintchine exponent $\psi$ satisfies $|\mbox{Im} (\psi)|\leq
M(1+\mbox{Re}(\psi))$ for some positive constant $M$.
Rao \cite{R77} gave a short proof of the Kanda-Forst theorem under the weaker condition that $X$ has resolvent densities.
 In particular, for $n\ge 1$, all
stable processes with index $\alpha\neq 1$ satisfy (H). Kanda
\cite{Ka78} proved that (H) holds for stable processes on
$\mathbf{R}^n$ with index $\alpha= 1$ if we assume that the linear
term vanishes. Silverstein {\cite{Si77} extended the Kanda-Forst
condition to the non-symmetric Dirichlet forms setting,
Fitzsimmons \cite{Fi01} extended it to the semi-Dirichlet forms
setting and Han et al. \cite{HMS11} extended it to the
positivity-preserving forms setting. Glover and Rao \cite{GR86}
proved that $\alpha$-subordinates of general Hunt processes
satisfy (H) (cf. Theorem \ref{GR} below).
{Rao \cite{R88} proved that if all 1-excessive functions of $X$ are
lower semicontinuous and $|{\rm Im}(\psi)|\leq (1+{\rm
Re}(\psi))f(1+{\rm Re}(\psi))$, where $f$ is an increasing function
on $[1,\infty)$ such that $\int_N^{\infty}(\lambda f(\lambda))^{-1}d\lambda=\infty$ for
every $N\geq 1$, then $X$ satisfies (H).}

Let $X$ be a L\'{e}vy process  on $\mathbf{R}^n$ with
L\'{e}vy-Khintchine exponent $(a,Q,\mu)$.  In \cite{HS11}, we showed that if $Q$ is non-degenerate then $X$ satisfies (H); if $Q$ is degenerate then, under the assumption that $\mu({\mathbf{R}^n\backslash
\sqrt{Q}\mathbf{R}^n})<\infty$, $X$ satisfies (H) if
and only if the equation
$$
\sqrt{Q}y=-a-\int_{\{x\in {\mathbf{R}^n\backslash
\sqrt{Q}\mathbf{R}^n}:\,|x|<1\}}x\mu(dx)
$$
has at least one solution $y\in \mathbf{R}^n$. We also showed that if $X$ is a
subordinator and satisfies (H) then its drift coefficient must be 0.

In this paper, we  will continue to explore (H) for L\'{e}vy
processes. The rest of the paper is organized as follows. In
Section 2, we present a comparison result on L\'{e}vy processes
which shows that big jumps have no effect on the validity of (H)
in some sense. Based on this result and the Kanda-Forst-Rao
theorem, in Section 3, we give examples of subordinators
satisfying (H). In Section 4, we give a new necessary and
sufficient condition for (H) and obtain an extended
Kanda-Forst-Rao theorem. By virtue of this theorem, we give a new
class of L\'{e}vy processes satisfying (H). In section 5, we
construct a type of subordinators that does not satisfy Rao's
condition. To the best of our knowledge, no existing criteria can
be applied to this example. It suggests that maybe new ideas and
methods are needed in order to completely solve Getoor's
conjecture even for the case of subordinators.

\section{A comparison result on  L\'{e}vy processes}\setcounter{equation}{0}
In this section, we prove a comparison result on L\'{e}vy
processes which implies that big jumps have no effect on the
validity of (H) in some sense.

Let $X$ be a L\'{e}vy process  on $\mathbf{R}^n$ with
L\'{e}vy-Khintchine exponent  $(a,Q,\mu)$. Suppose  that $\mu_1$
is a finite measure on $\mathbf{R}^n\backslash\{0\}$ such that
$\mu_1\leq \mu$. Denote $\mu_2:=\mu-\mu_1$ and let $X'$ be a
L\'{e}vy process on $\mathbf{R}^n$ with L\'{e}vy-Khintchine
exponent $(a',Q,\mu_2)$, where
\begin{eqnarray*}
a':=a+\int_{\{|x|<1\}}x\mu_1(dx).
\end{eqnarray*}

\begin{thm}\label{thm2.1}
Let $X$ and $X'$ be L\'{e}vy processes defined as above.
Then

\noindent (i) they have same semipolar sets.

\noindent (ii) they have same essentially polar sets.

\noindent (iii) if both $X$ and $X'$ have resolvent densities, then $X$
satisfies (H) if and only if  $X'$ satisfies (H).
\end{thm}
{\bf Proof.}  Denote by $\psi$ and $\psi'$ the L\'{e}vy-Khintchine
exponents of $X$ and $X'$, respectively. Then,
\begin{eqnarray}\label{thm2.1-a}
\psi'(z)&=&i\langle a',z\rangle+\frac{1}{2}\langle
z,Qz\rangle+\int_{\mathbf{R}^n} \left(1-e^{i\langle
z,x\rangle}+i\langle z,x\rangle
1_{\{|x|<1\}}\right)\mu_2(dx),\nonumber\\
\psi(z)&=&i\langle a,z\rangle+\frac{1}{2}\langle
z,Qz\rangle+\int_{\mathbf{R}^n} \left(1-e^{i\langle
z,x\rangle}+i\langle z,x\rangle
1_{\{|x|<1\}}\right)\mu(dx)\nonumber\\
&=&\psi'(z)+\int_{\mathbf{R}^n}\left(1-e^{i\langle
z,x\rangle}\right)\mu_1(dx).
\end{eqnarray}

(i) Suppose that $Y$ is a compound Poisson process with L\'{e}vy
measure $\mu_1$ and is independent of $X'$. By (\ref{thm2.1-a}),
$X$ has the same law as that of $X'+Y$. Let $T_1$ be the first
jumping time of $Y$. Then $T_1$ possesses an exponential
distribution and thus $P(T_1>0)=1$. Hence, for any set $A$ and any
point $x\in \mathbf{R}^n$, $x$ is a regular point of $A$ relative
to $X$ if and only if it is a regular point of $A$ relative to
$X'$. Therefore $X$ and $X'$ have same semipolar sets.

(ii) Set $C:=\mu_1(\mathbf{R}^n\backslash\{0\})$. By
(\ref{thm2.1-a}), we get
\begin{eqnarray}\label{sec2-3}
\mbox{Re}\psi'(z)\leq \mbox{Re}\psi(z)\leq
\mbox{Re}\psi'(z)+C
\end{eqnarray}
and
\begin{eqnarray}\label{sec2-4}
|\mbox{Im}\psi(z)|\leq |\mbox{Im}\psi'(z)|+C,\ \
|\mbox{Im}\psi'(z)|\leq |\mbox{Im}\psi(z)|+C.
\end{eqnarray}

For $\lambda>0$, we have
\begin{eqnarray}
\mbox{Re}\left(\frac{1}{\lambda+\psi(z)}\right)
&=&\frac{\lambda+\mbox{Re}\psi(z)}{(\lambda+\mbox{Re}\psi(z))^2+(\mbox{Im}\psi(z))^2},\label{sec2-5}\\
\mbox{Re}\left(\frac{1}{\lambda+\psi'(z)}\right)
&=&\frac{\lambda+\mbox{Re}\psi'(z)}{(\lambda+\mbox{Re}\psi'(z))^2+
(\mbox{Im}\psi'(z))^2}.\label{sec2-6}
\end{eqnarray}
By (\ref{sec2-3}) and (\ref{sec2-4}), we find that if $\lambda\geq
\sqrt{2}C$ then
\begin{eqnarray}\label{sec2-7}
\frac{\lambda+\mbox{Re}\psi(z)}{(\lambda+\mbox{Re}\psi(z))^2+(\mbox{Im}\psi(z))^2}
&\geq&
\frac{\lambda+\mbox{Re}\psi'(z)}{(\lambda+\mbox{Re}\psi'(z)+C)^2+
(|\mbox{Im}\psi'(z)|+C)^2}\nonumber\\
&\geq&\frac{\lambda+\mbox{Re}\psi'(z)}
{2[(\lambda+\mbox{Re}\psi'(z))^2+2C^2+(\mbox{Im}\psi'(z)^2]}\nonumber\\
&\geq&\frac{1}{4}\frac{\lambda+\mbox{Re}\psi'(z)}
{(\lambda+\mbox{Re}\psi'(z))^2+(\mbox{Im}\psi'(z))^2}.
\end{eqnarray}
Similar to (\ref{sec2-7}), we find that if $\lambda\geq 2C$ then
\begin{eqnarray}\label{sec2-8}
\frac{\lambda+\mbox{Re}\psi'(z)}{(\lambda+\mbox{Re}\psi'(z))^2+(\mbox{Im}\psi'(z))^2}
&\geq&
\frac{\lambda+\mbox{Re}\psi(z)-C}{(\lambda+\mbox{Re}\psi(z))^2+
(|\mbox{Im}\psi(z)|+C)^2}\nonumber\\
&\geq&\frac{\frac{1}{2}\lambda+\mbox{Re}\psi(z)}{(\lambda+\mbox{Re}\psi(z))^2+2C^2+
2(\mbox{Im}\psi(z))^2}\nonumber\\
&\geq&\frac{1}{4}\frac{\lambda+\mbox{Re}\psi(z)}{(\lambda+\mbox{Re}\psi(z))^2+(\mbox{Im}\psi(z))^2}.
\end{eqnarray}
By  (\ref{sec2-5})-(\ref{sec2-8}),  we obtain that if $\lambda\geq
2C$ then for any $z\in \mathbf{R}^n$,
\begin{eqnarray}\label{sec2-9}
\frac{1}{4}\ \mbox{Re}\left(\frac{1}{\lambda+\psi'(z)}\right)\leq
\mbox{Re}\left(\frac{1}{\lambda+\psi(z)}\right)\leq 4\
\mbox{Re}\left(\frac{1}{\lambda+\psi'(z)}\right).
\end{eqnarray}
By (\ref{sec2-9}) and Hawkes \cite[Theorem 3.3]{Ha79}, we obtain that a
set is essentially polar for $X$ if and only if it is essentially
polar for $X'$.

(iii) This is a direct consequence of  (i), (ii) and \cite[Theorem 2.1]{Ha79}.\hfill\fbox

For $\delta>0$, we define
$$
B_{\delta}:=\{x\in\mathbf{R}^n:0<|x|<
\delta\}.
$$

\begin{cor}\label{cor2.2}
Let $X_{\delta}$ be a L\'{e}vy process on $\mathbf{R}^n$
with L\'{e}vy-Khintchine exponent
$(a_{\delta},Q,\mu|_{B_{\delta}})$, where
\begin{eqnarray*}
a_{\delta}:=\left\{
\begin{array}{ll}
a+\int_{\{\delta\leq |x|<1\}}x\mu(dx),&\ \mbox{if}\ \ 0<\delta<1,\\
a,&\ \mbox{if}\ \ \delta\geq 1.
\end{array}
\right.
\end{eqnarray*}
Then, all the assertions of Theorem \ref{thm2.1} hold with $X'$ replaced by
$X_{\delta}$.
\end{cor}

\begin{rem}\label{rem2.3}
If  $\int_{|x|\leq 1}|x|\mu(dx)<\infty$, then $\psi$ can be
expressed by
$$\psi(z)=i\langle
d,z\rangle+\frac{1}{2}\langle z,Qz\rangle+\int_{\mathbf{R}^n}
\left(1-e^{i\langle z,x\rangle}\right)\mu(dx),
$$
where $-d$ is called the drift of $X$. In this case, we  call $(d,
Q,\mu)$ the L\'{e}vy-Khintchine exponent of $X$. For $\delta>0$,
we define $B_{\delta}$ and $X_{\delta}$ as above. Let
$X_{\delta}'$ be a L\'{e}vy process on $\mathbf{R}^n$ with
L\'{e}vy-Khintchine exponent $(d, Q,\mu|_{B_{\delta}})$. We claim
that $X_{\delta}$ and $X_{\delta}'$ have the same law and then all
the assertions of Theorem \ref{thm2.1} hold with $X'$ replaced by
$X_{\delta}'$. In fact, we have
\begin{eqnarray}\label{rem2.2-a}
d=a+\int_{\{|x|<1\}}x\mu(dx).
\end{eqnarray}
If $0<\delta<1$, then
\begin{eqnarray}\label{rem2.2-b}
a_{\delta}+\int_{\{|x|<1\}}x\mu|_{B_{\delta}}(dx)=\left(a+\int_{\{\delta\leq
|x|<1\}}\mu(dx)\right)+\int_{\{|x|<\delta\}}x\mu(dx)=d;
\end{eqnarray}
if $\delta\geq 1$, then
\begin{eqnarray}\label{rem2.2-c}
a_{\delta}+\int_{\{|x|<1\}}x\mu|_{B_{\delta}}(dx)=a+\int_{\{|x|<1\}}x\mu(dx)=d.
\end{eqnarray}
By (\ref{rem2.2-a})-(\ref{rem2.2-c}), we know that $X_{\delta}$
and $X_{\delta}'$ have the same L\'{e}vy-Khintchine exponent $(d,
Q,\mu|_{B_{\delta}})$ and thus have the same law.
\end{rem}


\section{Examples of subordinators satisfying (H)}\setcounter{equation}{0}
In this section, we will present new examples of subordinators
satisfying (H) by virtue of the comparison result given in Section
2 and the Kanda-Forst-Rao theorem. To the best of our knowledge,
which subordinators satisfy (H) is unknown in general. To
appreciate the importance of  the validity of (H) for
subordinators, let us recall the following remarkable result of
Glover and Rao.

\begin{thm}\label{GR} {\rm (Glover and Rao \cite{GR86})} Let $(X_t)_{t\ge 0}$ be a standard process on a locally compact space with a
countable base and $(T_t)_{t\ge 0}$ be an independent subordinator
satisfying Hunt's hypothesis (H). Then $(X_{T_t})_{t\ge 0}$
satisfies (H).
\end{thm}

Let $X$ be a subordinator. Then, its L\'{e}vy-Khintchine exponent
$\psi$ can be expressed by
$$\label{1-a}
\psi(z)=-idz+\int_{(0,\infty)}{\left(1-e^{izx}\right)}\mu(dx),\
z\in \mathbf{R},
$$
where $d\geq 0$ (called the drift coefficient) and $\mu$ satisfies
$\int_{(0,\infty)}(1\wedge x)\mu(dx)<\infty$. In \cite{HS11}, we
have proved the following result.

\begin{pro}\label{pro4.1}
If $X$ is a subordinator and satisfies (H), then $d=0$.
\end{pro}

By Proposition \ref{pro4.1}, when we consider (H) for
subordinators, we may concentrate on the case that $d=0$.
Hereafter we use $c_1, c_2,\dots$ to denote constants whose values
can change from one appearance to another.

\subsection{Special subordinators}
Let $X$ be a subordinator. Recall that the potential measure $U$
of $X$ is defined by
 $$
U(A)=E\left[\int_0^{\infty}1_{\{X_t\in A\}}dt\right],\ A\subset
[0,\infty).
 $$
For $\alpha>0$, the $\alpha$-potential measure $U^{\alpha}$ of $X$
is defined by
$$
U^{\alpha}(A)=E\left[\int_0^{\infty}e^{-\alpha t}1_{\{X_t\in
A\}}dt\right],\ A\subset [0,\infty).
$$

 $X$ is called a {\it special subordinator}  if  $U|_{(0,\infty)}$
has a decreasing density with respect to the Lebesgue measure.

\begin{thm}\label{thm4.2}
Let $X$ be a special subordinator. Then $X$ satisfies (H) if and
only if $d=0$.
\end{thm}
{\bf Proof.} By Proposition \ref{pro4.1}, we  need only prove the
sufficiency. Suppose that $d=0$. If $\mu$ is a finite measure,
then $X$ is a compound Poisson process and thus satisfies (H).

Now we consider the case that $\mu$ is an infinite measure. By
Bretagnolle \cite[Theorem 8]{Br71}, $X$ does not hit points, i.e., any single
point set $\{x\}$ is a polar set of $X$, which together with the
assumption that   $U|_{(0,\infty)}$ has a decreasing density with
respect to the Lebesgue measure, implies that $U|_{[0,\infty)}$
has a density with respect to the Lebesgue measure. Since for any
$\alpha>0$, $U^{\alpha}(\cdot)\leq U(\cdot)$, we obtain that for
any $\alpha\geq 0$, $U^{\alpha}$ is absolutely continuous with
respect to the Lebesgue measure. Then by Hawkes \cite[theorem 2.1]{Ha79},
we know that for any $\alpha\geq 0$, all $\alpha$-excessive
functions are lower semicontinuous. Therefore, by the fact that
$X$ does not hit points and  Blumenthal and Getoor \cite[Proposition (5.1), Theorem
(5.3)]{BG70}, following the same argument for stable subordinators
\cite[page 140]{BG70}, we obtain that $X$ satisfies
(H).\hfill\fbox

\subsection{Locally quasi-stable  subordinators}

  Let $S$ be a stable subordinator of index
$\alpha$, $0<\alpha<1$. Then, its L\'{e}vy-Khintchine exponent
$\psi_S$ has the form
$$\label{3.1}
\psi_S(z)=c|z|^{\alpha}(1-i\,\mbox{sgn}(z)\tan(\pi\alpha/2)),\
z\in (-\infty,\infty),
$$
where $c>0$. Its L\'{e}vy measure $\mu_S$ is absolutely continuous
with respect to the Lebesgue measure $dx$ and can be expressed by
\begin{eqnarray}\label{add001}
\mu_S(dx)=\left\{ \begin{array}{ll} c^+x^{-\alpha-1}dx,&\mbox{if}\
x>0,\\
0,&\mbox{if}\ x\le 0,
\end{array}
\right.
\end{eqnarray}
where $c^+>0$.

\begin{defi}\label{defi4.3}
Let $X$ be a subordinator with drift $0$ and L\'{e}vy measure
$\mu$. We call $X$ a {locally quasi-stable subordinator} if there
exist a stable subordinator $S$ with L\'{e}vy measure $\mu_S$,
positive constants $c_1,c_2,\delta$, and finite measures $\mu_1$
and $\mu_2$ on $(0,\delta)$ such that
$$c_1\mu_S-\mu_1\le
\mu\le c_2\mu_S+\mu_2\ \ {\rm on}\ (0,\delta). $$
\end{defi}

\begin{pro}\label{pro4.4}
Any locally quasi-stable subordinator satisfies (H).
\end{pro}
{\bf Proof.} Let $X,S,\mu_1$, $\mu_2$ and $\delta$ be as in
Definition \ref{defi4.3}. By Theorem \ref{thm2.1} and Remark
\ref{rem2.3}, we assume  without loss of generality that
$\mu|_{[\delta,\infty)}=0$ and  $\mu_1=0$. Denote by $\psi$ and
$\psi_S$ the L\'{e}vy-Khintchine exponents of $X$ and $S$,
respectively. Let $\mu_S$ be as in (\ref{add001}). Then
\begin{eqnarray}\label{pro5.4-a}
{\rm Re}\psi(z)&=&\int_0^{\infty}(1-\cos(zx))\mu(dx)\nonumber\\
&\ge&c_1\int_0^{\delta}(1-\cos(zx))\mu_S(dx)\nonumber\\
&=&c_1\left(\int_0^{\infty}(1-\cos(zx))\mu_S(dx)-\int_{\delta}^{\infty}(1-\cos(zx))\mu_S(dx)\right)\nonumber\\
 &=&c_1{\rm Re}\psi_S(z)-K_1\nonumber\\
 &=&c'|z|^{\alpha}-K_1,
\end{eqnarray}
where $c_1,c',K_1$ are positive constants.
\begin{eqnarray}\label{pro5.4-b}
|{\rm Im}\psi(z)|&\le&\int_0^{\infty}|\sin(zx)|\mu(dx)\nonumber\\
&\le&c_2\int_0^{\delta}|\sin(zx)|\mu_S(dx)\nonumber\\
&=&c_2\int_0^{\infty}|\sin(zx)|\mu_S(dx)-
c_2\int_{\delta}^{\infty}|\sin(zx)|\mu_S(dx)\nonumber\\
&\le&c_2c\left\{\int_0^{1/|z|}|\sin(zx)|x^{-1-\alpha}dx+
\int_{1/|z|}^{\infty}|\sin(zx)|x^{-1-\alpha}dx\right\}+K_2\nonumber\\
&\le&c_2c\left\{|z|\int_0^{1/|z|}x^{-\alpha}dx+\int_{1/|z|}^{\infty}x^{-1-\alpha}dx\right\}+K_2\nonumber\\
&=&c''|z|^{\alpha}+K_2,
\end{eqnarray}
where $c_2,c'',K_2$ are positive constants. By (\ref{pro5.4-a})
and (\ref{pro5.4-b}) we know that the Kanda-Forst condition holds
for $\psi$. By (\ref{pro5.4-a}) and Hartman and Wintner \cite{HW42}, we know that $X$
 has bounded continuous transition densities. Therefore, $X$ satisfies (H) by the Kanda-Forst theorem.\hfill\fbox

\begin{cor}\label{cor4.5}
Let $\varphi$ be a L\'{e}vy-Khintchine exponent and $\mu$ be a
L\'{e}vy measure of some special subordinator with drift 0 or some
locally quasi-stable subordinator. Then, the L\'{e}vy process with
L\'{e}vy-Khintchine exponent
\begin{eqnarray}\label{cor3.5-a}
\Phi(z):=\int_{(0,\infty)}\left(1-e^{-\varphi(z)x}\right) \mu(dx)
\end{eqnarray}
satisfies (H).
\end{cor}
{\bf Proof.} Let $X$ be a L\'{e}vy process with
L\'{e}vy-Khintchine exponent $\varphi$ and $(T_t)_{t\ge 0}$ be a
subordinator with drift 0 and L\'{e}vy measure $\mu$, which is
independent of $X$. Then $Y_t:=X_{T_t}$ has the L\'{e}vy exponent
$\Phi$ defined by (\ref{cor3.5-a}). Therefore, by Theorem
\ref{GR}, Theorem \ref{thm4.2} and Proposition \ref{pro4.4}, we
obtain that $Y$ satisfies (H).\hfill\fbox

\subsection{Further examples}
In this subsection, we give further examples of subordinators
satisfying  (H) by virtue of the comparison result given in
Section 2 and the following theorem of Rao.

\begin{thm}\label{R2} {\rm (Rao \cite{R88})} Let $X$ be a L\'evy process such that all 1-excessive functions are lower semicontinuous. Suppose there is an increasing function $f$ on $[1,\infty)$ such that
$ \int_N^{\infty}(\lambda f(\lambda))^{-1}d\lambda=\infty $ for
any $N\ge 1$ and  $|1+\psi|\le (1+{\rm Re}(\psi))f(1+{\rm
Re}(\psi))$. Then (H) holds.
\end{thm}

Let $0<\alpha<1$ and $0<\delta<1$. We define
$$
\mu_T(dx):=\frac{1}{-\log(x)x^{1+\alpha}}dx,\ \ 0<x<\delta
$$
and
$$
\mu_V(dx)=\frac{-\log(x)}{x^{1+\alpha}}dx,\ \ 0<x<\delta.
$$
 Let $X$ be a subordinator with drift 0 and L\'{e}vy measure
$\mu$.

\noindent (i) If $c_1\mu_T-\mu_1\le \mu\le c_2\mu_S+\mu_2$ on
$(0,\delta)$ for some positive constants $c_1,c_2$ and finite
measures $\mu_1,\mu_2$ on $(0,\delta)$, then $X$ satisfies (H).

In fact, by Theorem \ref{thm2.1} and Remark \ref{rem2.3}, we may
assume without loss of generality that $\mu|_{[\delta,\infty)}=0$
and $\mu_1=0$. For any $z\in \mathbf{R}$ with $|z|>1$, we have
\begin{eqnarray}\label{exm-1-a}
{\rm Re}\psi(z)&=&\int_0^{\infty}(1-\cos(zx))\mu(dx)\nonumber\\
&\ge&c_1\int_0^{\delta}(1-\cos(zx))\mu_T(dx)\nonumber\\
&=&c_1\int_0^{\infty}(1-\cos(zx))\mu_T(dx)-c_1\int_{\delta}^{\infty}(1-\cos(zx))\mu_T(dx)\nonumber\\
&\ge&c_1\int_{1/2|z|}^{1/|z|}(1-\cos(zx))\frac{1}{-\log(x)x^{1+\alpha}}dx-K_3\nonumber\\
&\ge&c_1'z^2\int_{1/2|z|}^{1/|z|}\frac{x^2}{-\log(x)x^{1+\alpha}}dx-K_3\nonumber\\
&\ge&c_1'\frac{z^2}{\log(2|z|)}\int_{1/2|z|}^{1/|z|}\frac{x^2}{x^{1+\alpha}}dx-K_3\nonumber\\
&=&c_1{''}\frac{|z|^{\alpha}}{\log(2|z|)}-K_3,
\end{eqnarray}
where $c_1',c_2'',K_3$ are positive constants. By (\ref{pro5.4-b})
and (\ref{exm-1-a}), we obtain that  $|{\rm Im}\psi(z)|\le
c^*(1+{\rm Re}\psi(z))\log(1+{\rm Re}\psi(z))$ for some positive
constant $c^*$. By Hartman and Wintner \cite{HW42} and (\ref{exm-1-a}), we know that
$X$
 has bounded continuous transition densities. Therefore, $X$ satisfies (H) by Theorem \ref{R2}.

 \bigskip

\noindent (ii) If $c_1\mu_S-\mu_1\le \mu\le c_2\mu_V+\mu_2$ on
$(0,\delta)$ for some positive constants $c_1,c_2$ and finite
measures $\mu_1,\mu_2$ on $(0,\delta)$, then $X$ satisfies (H).

In fact, by Theorem \ref{thm2.1} and Remark \ref{rem2.3}, we may
assume without loss of generality that $\mu|_{[\delta,\infty)}=0$
and $\mu_1=0$.  For any $z\in \mathbf{R}$ with $|z|>1/\delta$, we
have
\begin{eqnarray}\label{exm-2-a}
|{\rm
Im}\psi(z)|&\le&\int_0^{\infty}|\sin(zx)|\mu(dx)\nonumber\\
&\le& c_2\int_0^{\delta}|\sin(zx)|\mu_V(dx)
+K_4\nonumber\\
&\le&c_2\left\{\int_0^{1/|z|}|\sin(zx)|\frac{-\log(x)}{x^{1+\alpha}}dx+
\int_{1/|z|}^{\delta}|\sin(zx)|\frac{-\log(x)}{x^{1+\alpha}}dx\right\}+K_4\nonumber\\
&\le&c_2'\left\{|z|\int_0^{1/|z|}{\frac{-\log(x)}{x^{\alpha}}}dx+\log(|z|)
\int_{1/|z|}^{\infty}x^{-1-\alpha}dx\right\}+K_4\nonumber\\
&\le&c_2^{''}\left\{|z|\int_0^{1/|z|}{\frac{-(1-\alpha)\log(x)-1}
{x^{\alpha}}}dx+|z|^{\alpha}\log(|z|)\right\}+K_4\nonumber\\
&=&2c_2^{''}|z|^{\alpha}\log(|z|)+K_4,
\end{eqnarray}
where $c_2',c_2'',K_4$ are positive constants. By (\ref{pro5.4-a})
and (\ref{exm-2-a}), we obtain  that
 $|{\rm Im}\psi(z)|\le c^{**}{\rm Re}\psi(z)$ $\log({\rm
Re}\psi(z))$  for some positive constant $c^{**}$. By
(\ref{pro5.4-a}) and Hartman and Wintner \cite{HW42}, we know that $X$
 has bounded continuous transition densities. Therefore,  $X$ satisfies (H) by Theorem \ref{R2}.

\section{A new necessary and sufficient condition for (H) and an extended Kanda-Forst-Rao theorem}\setcounter{equation}{0}
Let $X$ be a L\'{e}vy process on $\mathbf{R}^n$. From now on we
assume that all 1-excessive functions are lower semicontinuous,
equivalently,  $X$ has resolvent densities. Define
$$
A:=1+{\rm Re}(\psi),\ \ B:=|1+\psi|.
$$

\begin{thm}\label{thm112} {\rm (Rao \cite{R88})} Let $\nu$ be a finite measure of finite 1-energy, i.e.,
$$\int_{\mathbf{R}^n} B^{-2}(z)A(z)|\hat \nu(z)|^2dz<\infty.$$ Then
\begin{equation}\label{H1}
\lim_{\lambda\rightarrow\infty}\int_{\mathbf{R}^n}|\hat\nu(z)|^2(\lambda+{\rm Re}\psi(z))|\lambda+\psi(z)|^{-2}dz
\end{equation}
exists. The limit is zero if and only if $U^1\nu$ is regular.
\end{thm}


%
Based on Theorems \ref{thm112} and \ref{R2}, we can prove the
following result.

\begin{lem}\label{Thm33} Let
$\nu$ be a finite measure of finite 1-energy and $f$ be an
increasing function on $[1,\infty)$ such that $
\int_N^{\infty}(\lambda f(\lambda))^{-1}d\lambda=\infty $ for some
$N\ge 1$. Then $U^1\nu$ is regular if and only if
$$
\lim_{\lambda\rightarrow\infty}\sum_{k=1}^{\infty}\int_{\{B(z)>A(z)f(A(z)),\,k\le \frac{|{\rm Im}\psi(z)|}{A(z)}<k+1,\,A(z)\le\lambda<(k+1)|{\rm Im}\psi(z)|\}}\frac{\lambda}{\lambda^2+({\rm Im}\psi(z))^2}|\hat \nu(z)|^2dz=0.
$$
\end{lem}

\noindent {\bf Proof.} Since $f$ is an increasing function on
$[1,\infty)$, $ \int_N^{\infty}(\lambda
f(\lambda))^{-1}d\lambda=\infty $ for some $N\ge 1$ if and only if
$ \int_N^{\infty}(\lambda f(\lambda))^{-1}d\lambda=\infty $ for
any $N\ge 1$. From the proof of Theorem \ref{R2} (see Rao \cite{R88}),
we know that the limit
\begin{equation}\label{gap2}
\lim_{\lambda\rightarrow\infty}\int_{A(z)\le\lambda}\frac{\lambda}{\lambda^2+B^2(z)}|\hat \nu(z)|^2dz
\end{equation}
exists and equals the limit in (\ref{H1}). We now show that the limit in (\ref{gap2}) equals 0 if and only if
\begin{equation}\label{H21}
\lim_{\lambda\rightarrow\infty}\int_{\{A(z)\le\lambda,\,B(z)>A(z)f(A(z))\}}\frac{\lambda}{\lambda^2+B^2(z)}|\hat \nu(z)|^2dz=0.
\end{equation}
To this end, we need only show that (\ref{H21}) implies that
\begin{equation}\label{gap}
\lim_{\lambda\rightarrow\infty}\int_{A(z)\le\lambda}\frac{\lambda}{\lambda^2+B^2(z)}|\hat \nu(z)|^2dz=0.
\end{equation}

Suppose that (\ref{H21}) holds. Then, the limit
$$
\lim_{\lambda\rightarrow\infty}\int_{\{A(z)\le\lambda,\,B(z)\le A(z)f(A(z))\}}\frac{\lambda}{\lambda^2+B^2(z)}|\hat \nu(z)|^2dz
$$
exists since the limit in (\ref{gap2}) always exists. Note that
\begin{eqnarray*}
& &\int_1^{\infty}\lambda^{-1}f(\lambda)^{-1}d\lambda\int_{\{A(z)\le\lambda,\,B(z)\le A(z)f(A(z))\}}\lambda(\lambda^2+B^2(z))^{-1}|\hat \nu(z)|^2dz\\
&=&\int_{\{B(z)\le A(z)f(A(z))\}}|\hat \nu(z)|^2dz\int_{A(z)}^{\infty}[f(\lambda)(\lambda^2+B^2(z))]^{-1}d\lambda\\
&\le&\frac{\pi}{2}\int_{\{B(z)\le A(z)f(A(z))\}}[B(z)f(A(z))]^{-1}|\hat \nu(z)|^2dz\\
&\le&\frac{\pi}{2}\int_{\mathbb{R}^d}B^{-2}(z)A(z)|\hat \nu(z)|^2dz\\
&<&\infty.
\end{eqnarray*}
Since $\int_1^{\infty}\lambda^{-1}f(\lambda)^{-1}d\lambda=\infty$,
$$\lim_{\lambda\rightarrow\infty}\int_{\{A(z)\le\lambda,\,B(z)\le A(z)f(A(z))\}}\frac{\lambda}{\lambda^2+B^2(z)}|\hat \nu(z)|^2dz=0.$$
Therefore, (\ref{gap}) holds by (\ref{H21}).

For each $k\in \mathbf{N}$, we have
\begin{eqnarray}\label{kljh1}
&&1_{\{k\le \frac{|{\rm Im}\psi(z)|}{A(z)}<k+1,\,\lambda\ge (k+1)|{\rm Im}\psi(z)|\}}\frac{\lambda}{\lambda^2+({\rm Im}\psi(z))^2}|\hat \nu(z)|^2\nonumber\\
&\le&1_{\{k\le \frac{|{\rm Im}\psi(z)|}{A(z)}<k+1,\,\lambda\ge (k+1)|{\rm Im}\psi(z)|\}}\frac{1}{\lambda}|\hat \nu(z)|^2\nonumber\\
&\le&\frac{1}{k+1}1_{\{k\le \frac{|{\rm Im}\psi(z)|}{A(z)}<k+1\}}\frac{|\hat \nu(z)|^2}{|{\rm Im}\psi(z)|}.
\end{eqnarray}
We assume without loss of generality that $f(1)=\sqrt{2}$. Note that $B(z)>A(z)f(A(z))$ implies that $B(z)\le\sqrt{2}|{\rm Im}\psi(z)|$. Then, we obtain by  $\int_{\mathbf{R}^n} B^{-2}(z)A(z)|\hat \nu(z)|^2dz<\infty$ that
\begin{equation}\label{kljh2}
\sum_{k=1}^{\infty}\frac{1}{2(k+1)}\int_{\{B(z)>A(z)f(A(z)),\,k\le \frac{|{\rm Im}\psi(z)|}{A(z)}<k+1\}} \frac{|\hat \nu(z)|^2}{|{\rm Im}\psi(z)|}dz<\infty.
\end{equation}
By (\ref{kljh1}), (\ref{kljh2}) and the dominated convergence theorem, we get
$$
\lim_{\lambda\rightarrow\infty}\sum_{k=1}^{\infty}\int_{\{B(z)>A(z)f(A(z)),\,k\le \frac{|{\rm Im}\psi(z)|}{A(z)}<k+1,\,\lambda\ge (k+1)|{\rm Im}\psi(z)|\}}\frac{\lambda}{\lambda^2+({\rm Im}\psi(z))^2}|\hat \nu(z)|^2dz=0.
$$
Therefore, the proof is complete by noting (\ref{H21}).\hfill\fbox

Note that if $\nu$ is a finite measure such that $U^1\nu$ is bounded then $\nu$ has finite 1-energy (cf.  Rao \cite[page 622]{R88}). By Lemma \ref{Thm33} and Proposition \ref{thm1}, we obtain the following necessary and sufficient condition for (H).

\begin{thm}\label{thmnb} Let $f$ be an increasing
function on $[1,\infty)$ such that $ \int_N^{\infty}(\lambda
f(\lambda))^{-1}d\lambda=\infty $ for some $N\ge 1$. Then (H)
holds if and only if
\begin{eqnarray}\label{thm2.4-1}
\lim_{\lambda\rightarrow\infty}\sum_{k=1}^{\infty}\int_{\{B(z)>A(z)f(A(z)),\,k\le \frac{|{\rm Im}\psi(z)|}{A(z)}<k+1,\,A(z)\le\lambda<(k+1)|{\rm Im}\psi(z)|\}}\frac{\lambda}{\lambda^2+({\rm Im}\psi(z))^2}|\hat \nu(z)|^2dz=0
\end{eqnarray}
for any finite measure $\nu$ with compact support such that $U^1\nu$ is bounded.
\end{thm}

\begin{rem}
Theorem \ref{thmnb} indicates that the validity of (H) is closely
related  to the behavior of $\psi(z)$ where ${\rm Im}(\psi(z))$ is
not well controlled by ${\rm Re}(\psi(z))$, which is possible and
can be seen from the uniform motion on $\mathbf{R}$ and the
example given in Section 5.\end{rem}

By virtue of Theorem \ref{thmnb}, we obtain the following result extending the Kanda-Forst-Rao theorem on (H).

\begin{thm}\label{cor11}
(H) holds if the following extended Kanda-Forst-Rao condition ((EKFR) for short) holds:\\
(EKFR) There are two measurable functions $\psi_1$ and $\psi_2$ on $\mathbf{R}^n$ such that $\rm{Im}(\psi)=\psi_1+\psi_2$,
and
\begin{eqnarray}\label{vbn1}
&&\quad\quad\quad |\psi_1|\leq Af(A),\nonumber\\
&&\int_{{\mathbf{R}^n}}\frac{|\psi_2(z)|}{(1+{\rm Re}\psi(z))^2+({\rm Im}\psi(z))^2}dz<\infty,
\end{eqnarray}
where  $f$ is an increasing function on $[1,\infty)$ such that $
\int_N^{\infty}(\lambda f(\lambda))^{-1}d\lambda=\infty $ for some
$N\ge 1$.
\end{thm}

\begin{rem}
If $\psi_2=0$, then the (EKFR) condition is just Rao's condition.
In particular, if $f=1$,  then it is just the Kanda-Forst
condition. An integrability condition similar to (\ref{vbn1}) has
been used in Glover \cite[Theorem 3.1]{Gl81}.
\end{rem}

\noindent {\bf Proof of Theorem \ref{cor11}.}
By Theorem \ref{thmnb}, we need only show that the limit in (\ref{thm2.4-1}) equals 0. We assume without loss of generality that $f(1)=1/3$. Note that $B(z)>3\sqrt{2}A(z)f(A(z))$ implies that $|{\rm Im}\psi(z)|> A(z)$ and $|{\rm Im}\psi(z)|> B(z)/\sqrt{2}$, and $|\psi_2(z)|>2A(z)f(A(z))$ implies that $|\psi_2(z)|>|{\rm Im}\psi(z)|/2$. Then,  by (\ref{vbn1}), the fact that $A(z)\leq c(1+|z|^2)$ for some positive constant $c$ and  the dominated convergence theorem, we obtain that
\begin{eqnarray*}
&&\sum_{k=1}^{\infty}\int_{\{B(z)>3\sqrt{2}A(z)f(A(z)),\,k\le \frac{|{\rm Im}\psi(z)|}{A(z)}<k+1,\,A(z)\le\lambda<(k+1)|{\rm Im}\psi(z)|\}}\frac{\lambda}{\lambda^2+({\rm Im}\psi(z))^2}|\hat \nu(z)|^2dz\\
&\le&\sum_{k=1}^{\infty}\int_{\{|{\rm Im}\psi(z)|>3A(z)f(A(z)),\,k\le \frac{|{\rm Im}\psi(z)|}{A(z)}<k+1,\,A(z)\le\lambda<(k+1)|{\rm Im}\psi(z)|\}}\frac{1}{2|{\rm Im}\psi(z)|}|\hat \nu(z)|^2dz\\
&\le&\sum_{k=1}^{\infty}\int_{\{|\psi_2(z)|>2A(z)f(A(z)),\,k\le \frac{|{\rm Im}\psi(z)|}{A(z)}<k+1,\,A(z)\le\lambda<(k+1)|{\rm Im}\psi(z)|\}}\frac{|\psi_2(z)|}{|{\rm Im}\psi(z)|^2}| \hat \nu(z)|^2dz\\
&\le&\sum_{k=1}^{\infty}\int_{\{k\le \frac{|{\rm Im}\psi(z)|}{A(z)}<k+1,\,A(z)\le\lambda<(k+1)|{\rm Im}\psi(z)|\}}\frac{2|\psi_2(z)|}{B^2(z)}|\hat \nu(z)|^2dz\\
&\le&\sum_{k=1}^{\infty}\int_{\{k\le \frac{|{\rm Im}\psi(z)|}{A(z)}<k+1,\,\lambda<(k+1)^2A(z)\}}\frac{2|\psi_2(z)|}{B^2(z)}|\hat \nu(z)|^2dz\\
&\le&\sum_{k=1}^{\infty}\int_{\{k\le \frac{|{\rm Im}\psi(z)|}{A(z)}<k+1,\,\lambda<c(k+1)^2(1+|z|^2)\}}\frac{2|\psi_2(z)|}{B^2(z)}|\hat \nu(z)|^2dz\\
&\rightarrow&0\ \ {\rm as}\ \lambda\rightarrow\infty.
\end{eqnarray*}
The proof is complete.\hfill\fbox

In the following, we give an application of Theorem \ref{cor11}.

\begin{thm}\label{new111} Let $\gamma>0$ and $X$ be
a L\'evy process on $ \mathbf{R}$ satisfying
\begin{equation}\label{new090}
\liminf_{|z|\rightarrow\infty}\frac{{\rm
Re}\psi(z)}{|z|\log^{\gamma}(|z|)}>0.
\end{equation}
Then $X$ satisfies (H).
\end{thm}
{\bf Proof.} By (\ref{new090}), we get
$$
\lim_{|z|\to\infty}\frac{{\rm Re}\psi(z)}{\log(1+|z|)}=\infty.
$$
Hence $X$ has bounded continuous transition densities by
Hartman and Wintner \cite{HW42}. Let $f(\lambda)=\log(\lambda)$ for $\lambda\in
[1,\infty)$ and set $\psi_1(z):=1_{\{|{\rm Im}\psi(z)|\le
A(z)f(A(z))\}}{\rm Im}\psi(z)$, $\psi_2(z):=1_{\{|{\rm
Im}\psi(z)|> A(z)f(A(z))\}}{\rm Im}\psi(z)$ for $z\in \mathbf{R}$.
Condition (\ref{new090}) implies that there exists a constant
$c>0$ such that
$$|\psi_2(z)|\ge c1_{\{|\psi(z)|>
A(z)f(A(z))\}}|z|\log^{1+\gamma}(|z|)$$ when $|z|$ is sufficiently
large. Therefore, (\ref{vbn1}) holds and the proof is complete by
Theorem \ref{cor11}.\hfill\fbox

\begin{exa} By Theorem \ref{new111}, Theorem \ref{thm2.1} and Corollary \ref{cor2.2}, we obtain a new class of 1-dimensional L\'{e}vy processes satisfying (H). Let $X$ be a L\'{e}vy process  on $\mathbf{R}$ with
L\'{e}vy-Khintchine exponent $(a,Q,\mu)$. Suppose that there exist
constants $\gamma>0$, $0<\delta<1$, $c>0$, and a finite measure
$\mu'$ on $\{x\in \mathbf{R}^n: 0<|x|<\delta\}$ such that
$$
d\mu\geq \frac{c(-\log(|x|))^{\gamma}}{x^2}dx-d\mu' \ \ {\rm  on}\
\{x\in \mathbf{R}: 0<|x|<\delta\}.$$ Similar to (\ref{exm-1-a}),
we can show that (\ref{new090}) holds. Then,  $X$ satisfies (H).
Note that in this example it does not matter if $a$ or $Q$ equals
0.

Let $Y$ be another 1-dimensional L\'{e}vy process which is
independent of $X$. Theorem \ref{new111} implies that the
perturbed process $Y+X$ also satisfies  (H).
\end{exa}

\begin{rem}
Blumenthal and Getoor introduced in \cite{BG61} the following
index $\beta''$ defined by
\begin{equation}\label{mnb}
\beta''=\sup\left\{\tau\geq 0:\frac{{\rm
Re}\psi(z)}{|z|^{\tau}}\to \infty\ \mbox{as}\ |z|\to
\infty\right\}.
\end{equation}
Let $X$ be a L\'evy process on $ \mathbf{R}$. Then, Theorem
\ref{new111} implies that (H) holds when $\beta''>1$. This result
is also a direct consequence of the following proposition.

\begin{pro}\label{pro-1} Let $X$ be
a L\'evy process on $ \mathbf{R}$. Suppose that
\begin{equation}\label{zxcv}
\liminf_{|z|\to \infty}\frac{|\psi(z)|}{|z|\log^{1+\gamma}|z|}>0
\end{equation}
for some constant $\gamma>0$. Then (H) holds.
\end{pro}
{\bf Proof.} Let $f\equiv 1$ and set $\psi_1(z):=1_{\{|{\rm
Im}\psi(z)|\le A(z)f(A(z))\}}{\rm Im}\psi(z)$,
$\psi_2(z):=1_{\{|{\rm Im}\psi(z)|> A(z)f(A(z))\}}{\rm Im}\psi(z)$
for $z\in \mathbf{R}$. Condition (\ref{zxcv}) implies that
\begin{eqnarray*}
\limsup_{|z|\to \infty}\left\{\frac{|\psi_2(z)|}{(1+{\rm
Re}\psi(z))^2+({\rm
Im}\psi(z))^2}\cdot|z|\log^{1+\gamma}|z|\right\}<\infty.
\end{eqnarray*}
Therefore, (\ref{vbn1}) holds and the proof is complete by Theorem
\ref{cor11}.\hfill\fbox

We remark that Proposition \ref{pro-1} can also be proved by
Theorem \ref{thm112}. In fact, the limit in (\ref{H1}) equals the
limit in (\ref{gap2}) and hence equals 0 by (\ref{zxcv}) and the
dominated convergence theorem.

\end{rem}

\section{A type of subordinators that does not satisfy Rao's condition}\setcounter{equation}{0}
As pointed out in Rao \cite{R88}, from the proof of Theorem \ref{R2}
it seems that the condition $B\le Af(A)$  is not far from being
necessary. In this section, however, we will construct a type of
subordinators that does not satisfy Rao's condition.

\subsection{Construction of the example}

We fix an $\alpha$ such that $\frac{1}{2}<\alpha<1$. In the
sequel, we define a function $\rho$ on $\mathbf{R}$ which will be
used as the density function of a L\'evy measure $\mu$.

First, we set $n_1=2$. Define a function $\rho_1$ on $\mathbf{R}$
as follows.
$$
\rho_1(x)=\frac{1}{x^{1+\alpha}},\ \ {\rm if}\ \frac{1}{2n^2_1}<x<\frac{1}{n^2_1};\ \ 0,\ {\rm otherwise}.
$$
We define $\mu_1(dx)=\rho_1(x)dx$ and denote by $\psi_1$ the
L\'evy-Khintchine exponent of $\mu_1$. Then, for $z\in
[\frac{n_1}{2}, 2n_1]$, we have
\begin{eqnarray}\label{jkl1}
{\rm Re}\psi_1(z)&=&\int_0^1(1-\cos(zx))\mu_1(dx)\nonumber\\
&\le&\frac{1}{2}\int_{1/2n_1^2}^{1/n_1^2}z^2x^2\frac{1}{x^{1+\alpha}}dx\nonumber\\
&\le&\frac{2n_1^{2\alpha-2}}{2-\alpha}\nonumber\\
&\le&2
\end{eqnarray}
and
\begin{eqnarray}\label{jkl2}
{\rm Im}\psi_1(z)&=&\int_0^1\sin(zx)\mu_1(dx)\nonumber\\
&=&\int_{1/2n_1^2}^{1/n_1^2}\sin(zx)\mu_1(dx)\nonumber\\
&\ge&\int_{1/2n_1^2}^{1/n_1^2}\frac{zx}{2x^{1+\alpha}}dx\nonumber\\
&\ge&\frac{1}{8}n_1^{2\alpha-1}.
\end{eqnarray}
We increase $n_1$ so that
$\frac{1}{8}n_1^{2\alpha-1}>\frac{6}{1-\alpha}$.

For any $z\in \mathbf{R}$, we have
\begin{equation}\label{21}
{\rm Re}\psi_1(z)=\int_0^1(1-\cos(zx))\mu_1(dx)\le
\int_{1/2n_1^2}^1\frac{1}{x^{1+\alpha}}dx\le
\frac{2^{\alpha}n_1^{2\alpha}}{\alpha}\le 4n_1^{2\alpha}
\end{equation}
and
\begin{equation}\label{22}
|{\rm Im}\psi_1(z)|\le\int_0^1|\sin(zx)|\mu_1(dx)\le
\int_{1/2n_1^2}^1\frac{1}{x^{1+\alpha}}dx\le
\frac{2^{\alpha}n_1^{2\alpha}}{\alpha}\le 4n_1^{2\alpha}.
\end{equation}

\vskip 0.3cm We choose an $n_2\in \mathbf{N}$ such that
$n^2_2>2n_1^2$. We define a function $\rho_2$ on $\mathbf{R}$ as
follows.
$$
\rho_2(x)=\frac{1}{x^{1+\alpha}},\ \ \mbox{if}\ \frac{1}{2n^2_2}<x<\frac{1}{n^2_2};\ \ 0,\ {\rm otherwise}.
$$
Note that there is no overlap between $\rho_1$ and $\rho_2$. We
define $\mu_2(dx)=\rho_2(x)dx$ and denote by $\psi_2$  the
L\'evy-Khintchine  exponent of $\mu_2$. Then, similar to the
above, we can show that for $z\in [\frac{n_2}{2}, 2n_2]$
\begin{equation}\label{1}
{\rm Re}\psi_2(z)\le 2\ \ {\rm and}\ \ {\rm Im}\psi_2(z)\ge
\frac{1}{8}n_2^{2\alpha-1}\left(>\frac{6}{1-\alpha}\right).
\end{equation}

Note that for $z\in [\frac{n_1}{2}, 2n_1]$ we have
\begin{eqnarray}\label{p1}
{\rm Re}\psi_2(z)&=&\int_0^1(1-\cos(zx))\mu_2(dx)\nonumber\\
&\le&\frac{1}{2}\int_{1/2n_2^2}^{1/n_2^2}z^2x^2\frac{1}{x^{1+\alpha}}dx\nonumber\\
&\le&\frac{2n_1^2n_2^{2\alpha-4}}{2-\alpha}
\end{eqnarray}
and
\begin{eqnarray}\label{p2}
|{\rm Im}\psi_2(z)|&\le&\int_0^1|\sin(zx)|\mu_2(dx)\nonumber\\
&\le&\int_{1/2n_2^2}^{1/n_2^2}|\sin(zx)|\frac{1}{x^{1+\alpha}}dx\nonumber\\
&\le&\int_{1/2n_2^2}^{1/n_2^2}2n_1x\frac{1}{x^{1+\alpha}}dx\nonumber\\
&\le&\frac{2n_1n_2^{2\alpha-2}}{1-\alpha}.
\end{eqnarray}
We increase $n_2$ (with $n_1$ fixed) so that $n_2\ge
n_1^{5/(2-2\alpha)}$. By (\ref{p1}) and (\ref{p2}), we get
\begin{equation}\label{az1}
{\rm Re}\psi_2(z)\le\frac{2}{(1-\alpha)n_1^4},\ \ |{\rm
Im}\psi_2(z)|\le\frac{2}{(1-\alpha)n_1^4},\ \ z\in
\left[\frac{n_1}{2}, 2n_1\right].
\end{equation}
Then, by (\ref{jkl1}), (\ref{jkl2}) and (\ref{az1}), we obtain
that for $z\in [\frac{n_1}{2}, 2n_1]$,
\begin{equation}\label{r1}
{\rm Re}\psi_1(z)+{\rm Re}\psi_2(z)\le 2+\frac{2}{(1-\alpha)n_1^4}
\end{equation}
and
\begin{equation}\label{r2}
{\rm Im}\psi_1(z)+{\rm Im}\psi_2(z)\ge
\frac{1}{8}n_1^{2\alpha-1}-\frac{2}{(1-\alpha)n_1^4}.
\end{equation}

We further increase $n_2$ so that $n_2\ge
(96)^{1/(2\alpha-1)}n_1^{(4+2\alpha)/(2\alpha-1)}$ which ensures
that  for any $z\in \mathbf{R}$ (cf. (\ref{21}), (\ref{22}) and
(\ref{1})),
\begin{equation}\label{az2}
{\rm Re}\psi_1(z)\le \frac{1}{3n_1^4}{\rm
Im}\psi_2\left(\frac{n_2}{2}\right),\ \ |{\rm Im}\psi_1(z)|\le
\frac{1}{3n_1^4}{\rm Im}\psi_2\left(\frac{n_2}{2}\right).
\end{equation}
By (\ref{1}) and (\ref{az2}), we obtain that for $z\in
[\frac{n_2}{2}, 2n_2]$,
\begin{equation}\label{q1}
{\rm Re}\psi_1(z)+{\rm Re}\psi_2(z)\le \frac{1}{3n_1^4}{\rm
Im}\psi_2\left(\frac{n_2}{2}\right)+2
\end{equation}
and
\begin{equation}\label{q2}
{\rm Im}\psi_1(z)+{\rm Im}\psi_2(z)\ge
\left(1-\frac{1}{3n_1^4}\right){\rm
Im}\psi_2\left(\frac{n_2}{2}\right).
\end{equation}
Define
\begin{equation}\label{q20}\vartheta:=\max\left\{\frac{5}{2-2\alpha},
\frac{4+2\alpha}{2\alpha-1}\right\}.\end{equation} We can set
$n_2$ to be $cn_1^{\vartheta}$, for some positive constant $c$
depending only on $\alpha$, such that (\ref{r1}), (\ref{r2}),
(\ref{q1}) and (\ref{q2}) hold.

For any $z\in \mathbf{R}$, we have
\begin{equation}\label{31}
{\rm Re}\psi_2(z)=\int_0^1(1-\cos(zx))\mu_2(dx)\le
\int_{1/2n_2^2}^1\frac{1}{x^{1+\alpha}}dx\le
\frac{2^{\alpha}n_2^{2\alpha}}{\alpha}\le4n_2^{2\alpha}
\end{equation}
and
\begin{equation}\label{32}
|{\rm Im}\psi_2(z)|\le\int_0^1|\sin(zx)|\mu_2(dx)\le
\int_{1/2n_2^2}^1\frac{1}{x^{1+\alpha}}dx\le
\frac{2^{\alpha}n_2^{2\alpha}}{\alpha}\le4n_2^{2\alpha}.
\end{equation}
\vskip 0.3cm We choose an $n_3\in \mathbf{N}$  such that
$n^3_2>2n_2^2$. We define a function $\rho_3$ on $\mathbf{R}$ as
follows.
$$
\rho_3(x)=\frac{1}{x^{1+\alpha}},\ \ \mbox{if}\ \frac{1}{2n^2_3}<x<\frac{1}{n^2_3};\ \ 0,\ {\rm otherwise}.
$$
Note that there is no overlap among $\rho_1$, $\rho_2$ and
$\rho_3$. We define $\mu_3(dx)=\rho_3(x)dx$  and denote by
$\psi_3$  the L\'evy-Khintchine exponent of $\mu_3$. Then, similar
to the above, we can show that  for $z\in [\frac{n_3}{2}, 2n_3]$,
\begin{equation}\label{41}
{\rm Re}\psi_3(z)\le 2\ \ {\rm and}\ \ {\rm Im}\psi_3(z)\ge
\frac{1}{8}n_3^{2\alpha-1}
\end{equation}
and for any $z\in \mathbf{R}$,
$$
{\rm Re}\psi_3(z)\le4n_3^{2\alpha},\ \ |{\rm Im}\psi_3(z)|\le
4n_3^{2\alpha}.
$$

Similar to (\ref{p1}) and (\ref{p2}), we obtain that  for $z\in
[\frac{n_1}{2}, 2n_1]$,
\begin{equation}\label{v1}
{\rm Re}\psi_3(z)\le\frac{2n_1^2n_3^{2\alpha-4}}{2-\alpha},\ \
|{\rm Im}\psi_3(z)|\le\frac{2n_1n_3^{2\alpha-2}}{1-\alpha}
\end{equation}
and for $z\in [\frac{n_2}{2}, 2n_2]$,
\begin{equation}\label{v2}
{\rm Re}\psi_3(z)\le\frac{2n_2^2n_3^{2\alpha-4}}{2-\alpha},\ \
|{\rm Im}\psi_3(z)|\le\frac{2n_2n_3^{2\alpha-2}}{1-\alpha}.
\end{equation}
We increase $n_3$  (with $n_1,n_2$ fixed) so that $n_3\ge
n_2^{5/(2-2\alpha)}$. By (\ref{v1}) and (\ref{v2}), we get
\begin{equation}\label{az22}
{\rm Re}\psi_3(z)\le\frac{2}{(1-\alpha)n_2^4},\ \ |{\rm
Im}\psi_3(z)|\le\frac{2}{(1-\alpha)n_2^4},\ \ z\in
\left[\frac{n_1}{2}, 2n_1\right]\bigcup\left[\frac{n_2}{2},
2n_2\right].
\end{equation}
Hence, by (\ref{r1}), (\ref{r2}) and (\ref{az22}), we obtain that
for $z\in [\frac{n_1}{2}, 2n_1]$,
\begin{equation}\label{rr1}
{\rm Re}\psi_1(z)+{\rm Re}\psi_2(z)+{\rm Re}\psi_3(z)\le
2+\frac{2}{(1-\alpha)n_1^4}+\frac{2}{(1-\alpha)n_2^4}
\end{equation}
and
\begin{equation}\label{rr2}
{\rm Im}\psi_1(z)+{\rm Im}\psi_2(z)+{\rm Im}\psi_3(z)\ge
\frac{1}{8}n_1^{2\alpha-1}-\frac{2}{(1-\alpha)n_1^4}-\frac{2}{(1-\alpha)n_2^4}.
\end{equation}
By (\ref{q1}), (\ref{q2}), (\ref{az22}) and (\ref{1}), we obtain
that for $z\in [\frac{n_2}{2}, 2n_2]$,
\begin{equation}\label{qq1}
{\rm Re}\psi_1(z)+{\rm Re}\psi_2(z)+{\rm Re}\psi_3(z)\le
\frac{2}{3n_1^4}{\rm
Im}\psi_2\left(\frac{n_2}{2}\right)+2+\frac{2}{(1-\alpha)n_2^4}
\end{equation}
and
\begin{equation}\label{qq2}
{\rm Im}\psi_1(z)+{\rm Im}\psi_2(z)+{\rm Im}\psi_3(z)\ge
\left(1-\frac{1}{3n_1^4}-\frac{1}{3n_2^4}\right){\rm
Im}\psi_2\left(\frac{n_2}{2}\right).
\end{equation}

We further increase $n_3$ so that $n_3\ge
(192)^{1/(2\alpha-1)}n_2^{(4+2\alpha)/(2\alpha-1)}$ which ensures
that  for any $z\in \mathbf{R}$ (cf. (\ref{21}), (\ref{22}),
(\ref{31}), (\ref{32}) and (\ref{41})),
\begin{equation}\label{az2223}
{\rm Re}\psi_1(z),{\rm Re}\psi_2(z),|{\rm Im}\psi_1(z)|,|{\rm
Im}\psi_2(z)|\le \frac{1}{6n_2^4}{\rm
Im}\psi_2\left(\frac{n_3}{2}\right).
\end{equation}
Therefore, we obtain by (\ref{41}) and (\ref{az2223}) that for
$z\in [\frac{n_3}{2}, 2n_3]$,
\begin{equation}\label{cq1}
{\rm Re}\psi_1(z)+{\rm Re}\psi_2(z)+{\rm Re}\psi_3(z)\le
\frac{1}{3n_2^4}{\rm Im}\psi_3\left(\frac{n_3}{2}\right)+2
\end{equation}
and
\begin{equation}\label{cq2}
{\rm Im}\psi_1(z)+{\rm Im}\psi_2(z)+{\rm Im}\psi_3(z)\ge
\left(1-\frac{1}{3n_1^4}-\frac{1}{3n_2^4}\right){\rm
Im}\psi_3\left(\frac{n_3}{2}\right).
\end{equation}
We set $n_3$ to be $ 2^{1/(2\alpha-1)}cn_2^{\vartheta}$, where
$\vartheta$ and $c$ are as the same as above.

\vskip 0.3cm Continue in this way, we define $\rho_4,\rho_5,\dots$
All of these functions have no overlap and we have estimates
similar to (\ref{rr1})-(\ref{qq2}), (\ref{cq1}) and (\ref{cq2}).
Now we define
$$
\rho=\sum_{i=1}^{\infty}\rho_i.
$$
One finds that $\mu(dx)=\rho(x)dx$ is the L\'evy measure of a
subordinator $X$ with the L\'evy-Khintchine exponent
$$
\psi=\sum_{i=1}^{\infty}\psi_i.
$$
Moreover, we have that for $k\ge 2$,
\begin{equation}\label{hgf}
n_k=(k-1)^{1/(2\alpha-1)}cn_{k-1}^{\vartheta},
\end{equation} and for $z\in [\frac{n_k}{2},2n_k]$,
\begin{equation}\label{hgf7}
{\rm Im}\psi_k(z)\ge\frac{1}{8}n_k^{2\alpha-1},
\end{equation}
\begin{equation}\label{general1}
{\rm Re}\psi(z)\le \frac{1}{3n_{k-1}^4}{\rm
Im}\psi_k\left(\frac{n_k}{2}\right)+2+\frac{2}{1-\alpha}\sum_{k=1}^{\infty}\frac{1}{n_k^4},
\end{equation}
and
\begin{equation}\label{general2}
{\rm Im}\psi(z)\ge
\left(1-\frac{1}{3}\sum_{k=1}^{\infty}\frac{1}{n_k^4}\right){\rm
Im}\psi_k\left(\frac{n_k}{2}\right).
\end{equation}

\subsection{Discussions}

In this subsection, we make  discussion about  the subordinators
constructed in Subsection 5.1. Below we use $c_1, c_2,\dots$ to
denote positive constants depending only on $\alpha$.

\noindent \textbf{1.} By the estimates (\ref{general1}) and
(\ref{general2}), we can show that Rao's condition does not hold
for the subordinators. In fact, by (\ref{hgf}), there exists a
constant $c_1> 1$ such that
\begin{equation}\label{explain}
n_k>c_1^{c_1^k}, \ \ k\in \mathbf{N}. \end{equation}
 By (\ref{hgf7}),
(\ref{general1}) and (\ref{general2}), we find that there exist
constants $c_2,c_3,c_4>0$  such that for any $k\ge2$,
\begin{equation}\label{add345}
\frac{{\rm Im}\psi(z)}{1+{\rm Re}\psi(z)}\ge c_2n_{k-1}^4\ge
c_3n_k^{3/\vartheta}\ge
c_3\left(\frac{z}{2}\right)^{3/\vartheta},\ \ \forall z\in[n_k/2,
2n_k].
\end{equation}
\begin{equation}\label{break1}
{\rm Re}\psi(z)\le c_4n_{k-1}^{\alpha\vartheta-3},\ \ \forall
z\in[n_k/2, 2n_k].
\end{equation}
The estimates (\ref{add345}) and (\ref{break1}) imply that there
does not exist an increasing function $f$ on $[1,\infty)$
satisfying $ \int_N^{\infty}(\lambda
f(\lambda))^{-1}d\lambda=\infty $ for some $N\ge 1$ and
$|1+\psi|\le (1+{\rm Re}(\psi))f(1+{\rm Re}(\psi))$. That is,
Rao's condition does not hold for the subordinators constructed in
Subsection 5.1.

By Theorem \ref{thm2.1}, we can modify the L\'evy measure $\mu$
defined in Subsection 5.1 by a finite measure and hence obtain a
subordinator which does not satisfy Rao's condition and whose
L\'evy measure $\mu$ has a smooth density $\rho$ with respect to
the Lebesgue measure on $(0,\infty)$.

\vskip 0.3cm \noindent \textbf{2.} Besides the index $\beta''$
(see (\ref{mnb})), Blumenthal and Getoor introduced also in
\cite{BG61} the indexes $\beta$ and $\sigma$ defined by
$$
\beta=\inf\left\{\tau>0:
\int_{\{|x|<1\}}|x|^{\tau}\mu(dx)<\infty\right\}
$$
and
$$
\sigma=\sup\left\{\tau\le 1:
\int_1^{\infty}\frac{x^{\tau-1}}{\int_0^{\infty}(1-e^{-xy})\mu(dy)}dx<\infty\right\}.
$$

From the construction of the subordinators given in Subsection
5.1, we obtain by \cite[Theorem 6.1]{BG61} that
$$
\sigma=\beta=\alpha.
$$
By (\ref{hgf}) and (\ref{general1}) (cf.  (\ref{pro5.4-b})), we
get
$$
\beta''\le \alpha-\frac{4}{\vartheta}.
$$

\vskip 0.1cm \noindent \textbf{3.} Take $\alpha=3/4$. For the
subordinators constructed in Subsection 5.1, we claim that there
exists a finite signed measure $d\nu=g_1dx-g_2dx$ with $g_1,g_2\in
L^1_+(\mathbb{R};dx)$ such that
\begin{equation}\label{pop11}\int_{\mathbf{R}}
B^{-2}(z)A(z)|{\hat\nu}(z)|^2dz<\infty\end{equation}
 but
\begin{equation}\label{pop12}
\lim_{\lambda\rightarrow\infty}\int_{\mathbf{R}}|{\hat{\nu}}(z)|^2(\lambda+{\rm
Re}\psi(z))|\lambda+\psi(z)|^{-2}dz=\infty.\end{equation}

Let $\omega$ be a sufficiently large number. We define
$$
\zeta_{\omega}(x):=\left\{1-\frac{1-1/({\omega})^{0.1}}{{\omega}}\cdot|x|\right\},
\ \ {\rm if}\ |x|\le\omega;\ \ \frac{1}{|x|^{0.1}},\ {\rm
otherwise},
$$
and
$$
\eta_{\omega}(x):=\left\{1-\frac{1-1/({\omega})^{0.1}}{{\omega}}\cdot|x|\right\}\vee
0, \ \ x\in \mathbb{R}.
$$
By Polya's theorem (cf. Lukacs \cite[Theorem 4.3.1]{Luka}), both
$\zeta_{\omega}$ and $\eta_{\omega}$ are characteristic functions
of absolutely continuous symmetric distributions. Define
$\varsigma_{\omega}:=\eta_{\omega}-\zeta_{\omega}$. Then, $
\varsigma_{\omega}(x)=0$ if $|x|\le{\omega}$; $
\varsigma_{\omega}(x)=1/|x|^{0.1}$ if $|x|\ge (1.1){\omega}$; and
$0\le \varsigma_{\omega}(x)\le 1/|x|^{0.1}$ otherwise.

Let $k_0\in \mathbf{N}$ be a sufficiently large number. For $k\ge
k_0$, we define
$\xi_k:=\varsigma_{\frac{n_k}{2}}-\varsigma_{\frac{2n_k}{1.1}}$.
We find that $\xi_k$ is a characteristic function of the
difference of two functions $g^k_1,g^k_2\in L^1_+(\mathbb{R};dx)$
with $\|g_1^k\|_{L^1},\|g_2^k\|_{L^1}\le 2$. Define
$g_1:=\sum_{k=1}^{\infty}g^k_1/2^k$,
$g_2:=\sum_{k=1}^{\infty}g^k_2/2^k$ and $d\nu:=g_1dx-g_2dx$. By
applying (\ref{q20}), (\ref{hgf}), (\ref{explain}) and the first
inequality of (\ref{add345}) to $B(z)/A(z)$ and applying
(\ref{hgf7}), (\ref{general2}) to $B(z)$, we find that there
exists a constant $c_5>0$ such that
\begin{eqnarray*}
\int_{\mathbf{R}}
B^{-2}(z)A(z)|{\hat\nu}(z)|^2dz
&=&\int_{\mathbf{R}}
\frac{1}{\frac{B(z)}{A(z)}\cdot B(z)}|{\hat\nu}(z)|^2dz\\
&\le&c_{5}\sum_{k=1}^{\infty}\frac{1}{n_k^{\frac{4}{\vartheta}-\frac{1}{22}}\cdot
n_k^{2\alpha-1}\cdot
2^{2k}}\int_{n_k/2}^{2n_k}\frac{1}{z^{0.2}}dz\\
&=&c_{5}\sum_{k=1}^{\infty}\frac{1}{n_k^{9/11}\cdot
2^{2k}}\int_{n_k/2}^{2n_k}\frac{1}{z^{0.2}}dz\\
&<&\infty.\end{eqnarray*}
However, there exists a constant $c_6>0$
 such that (cf. ({\ref{pro5.4-b} and (\ref{explain}))
\begin{eqnarray*}
\int_{\mathbf{R}}
|{\hat\nu}(z)|^2\frac{n_k^{\alpha}}{(n_k^{\alpha})^2+({\rm
Im}\psi(z))^2}dz&\ge&c_{6}\frac{1}{n_k^{\frac{3}{4}}\cdot 2^{2k}}\int_{(0.55)n_k}^{\frac{2n_k}{1.1}}\frac{1}{z^{0.2}}dz\\
&\rightarrow&\infty\ \ {\rm as}\
k\rightarrow\infty,\end{eqnarray*} which implies (\ref{pop12}).

By (\ref{pop11}) and (\ref{pop12}) we can also conclude that Rao's
condition does not hold for the subordinators constructed in
Subsection 5.1. In fact, from the proof of Theorem \ref{R2} (see
Rao \cite{R88}), we can see that under Rao's condition, $$
\lim_{\lambda\rightarrow\infty}\int_{\mathbf{R}}|{\hat{\nu}}(z)|^2(\lambda+{\rm
Re}\psi(z))|\lambda+\psi(z)|^{-2}dz=0
$$
holds for any finite signed measure of finite 1-energy.

 It is interesting to
compare (\ref{pop11}) and (\ref{pop12}) with the following result,
which is a consequence of Theorem \ref{thm112}.
\begin{thm}\label{slight} Let $X$ be a L\'evy process on $\mathbf{R}^n$ such that all 1-excessive functions are lower
semicontinuous. Then (H) holds if and only if
\begin{equation}\label{sx}
\lim_{\lambda\rightarrow\infty}\int_{\mathbf{R}^n}|\hat\nu(z)|^2(\lambda+{\rm
Re}\psi(z))|\lambda+\psi(z)|^{-2}dz=0
\end{equation}
for any finite measure $\nu$ of finite 1-energy.
\end{thm}
{\bf Proof.} By Theorem \ref{thm112}, Rao \cite[Remark, page 622]{R88}
and Blumenthal and Getoor \cite[VI. (4.8)]{BG68}, we need only prove the necessity.
Suppose that (H) holds for $X$. Let $\nu$ be a finite measure of
finite 1-energy and $\kappa$ be the standard Gaussian measure on
$\mathbf{R}^n$. Then, $\nu+\kappa$ has finite 1-energy, which
implies that
\begin{equation}\label{pop13}
\int_{\mathbf{R}^n}U^1(\nu+\kappa)d(\nu+\kappa)<\infty.
\end{equation}
By (\ref{pop13}), $\kappa(\{x:U^1(\nu+\kappa)(x)=\infty\})=0$.
Hence $U^1(\nu+\kappa)$ is locally integrable (with respect to the
Lebesgue measure $dx$) by \cite[VI. (2.3)]{BG68}. By (H) and
\cite[VI. (4.9)]{BG68}, we find that $U^1(\nu+\kappa)$ is regular.
Therefore, (\ref{sx}) holds by Theorem \ref{thm112} and the proof
is complete.\hfill\fbox

 So far we have not been able to
prove or disprove that (H) holds for the subordinators constructed
in Subsection 5.1. This example suggests that maybe completely new
ideas and methods are needed for resolving Getoor's conjecture.

\bigskip

{ \noindent {\bf\large Acknowledgments} \vskip 0.1cm  \noindent We
thank Professor Fengyu Wang for helpful comments that improved a
previous version of the paper.  We are grateful to the support of
NNSFC, Jiangsu Province Basic Research Program (Natural Science
Foundation) (Grant No. BK2012720), and NSERC.}


\begin{thebibliography}{1234}



\bibitem{BG61} Blumenthal R.M., Getoor R.K.: Sample functions of stochastic processes with stationary independent increments. J. Math. Mech., \textbf{10}, 493-516 (1961).

\bibitem{BG68} Blumenthal R.M., Getoor R.K.:  Markov Processes and
Potential Theory. Academic Press, New York and London (1968).


\bibitem{BG70}Blumenthal R.M., Getoor R.K.:  Dual processes and
potential theory. Proc. 12th Biennial Seminar of the Canadian Math.
Congress, 137-156 (1970).

\bibitem{Br71} Bretagnolle J.:  R\'{e}sults de Kesten sur les processus \`{a}
accroissements ind\'{e}pendants. S\'{e}minare de Probabilit\'{e}s V,
Lect. Notes in Math., Vol. 191, Springer-Verlag, Berlin, 21-36
(1971).

\bibitem{Fi90}Fitzsimmons P.J.: On the equivalence of three
potential principles for right Markov processes. Probab. Th. Rel.
Fields \textbf{84}, 251-265 (1990).

\bibitem{Fi01} Fitzsimmons P.J.:  On the quasi-regularity of semi-Dirichlet
forms. Potential Anal. {\bf 15}, 151-185 (2001).

\bibitem{FK} Fitzsimmons P.J., Kanda M.:  On Choquet's dichotomy of capacity for Markov processes.
Ann. Probab. {\bf 20}, 342-349 (1992).

\bibitem{F75}Forst
G.:  The definition of energy in non-symmetric translation
invariant Dirichlet spaces. Math. Ann. {\bf 216}, 165-172 (1975).

%
%


\bibitem{Gl81} Glover J.: Energy and the maximum principle for
nonsymmetric Hunt processes. Probability Theory and Its
Applications, XXVI, 4, 757-768 (1981).


\bibitem{G83}Glover J.: Topics in energy and potential theory. Seminar on Stochastic Processes, 1982, Birkh\"{a}user, 195-202 (1983).

%
\bibitem{GR86} Glover J., Rao M.:  Hunt's hypothesis (H) and
Getoor's conjecture. Ann. Probab. {\bf 14}, 1085-1087 (1986).

\bibitem{HMS11}Han X.-F., Ma Z.-M. and Sun W.: $h\hat{h}$-transforms of positivity preserving semigroups and associated Makov processes.
Acta Math. Sinica, English Series \textbf{27}, 369-376 (2011).


\bibitem{HW42}  Hartman P., Wintner A.:   On the infinitesimal
generators of integral convolutions. Amer. J. Math. {\bf 64}, 273-298 (1942).



\bibitem{Ha79} Hawkes J.: Potential theory of L\'{e}vy processes. Proc. London Math. Soc. {\bf 3},
335-352 (1979).

\bibitem{HS11} Hu Z.-C. and Sun W.: Hunt's hypothesis (H) and
Getoor's conjecture for L\'{e}vy processes. Stoch. Proc. Appl.
\textbf{122}, 2319-2328 (2012).

%
%
\bibitem{Ka76} Kanda M.:  Two theorems on capacity for Markov processes
with stationary independent increments. Z. Wahrsch. verw. Gebiete {\bf 35}, 159-165 (1976).
%
\bibitem{Ka78} Kanda M.:  Characterisation of semipolar sets
for processes with stationary independent increments. Z. Wahrsch.
verw. Gebiete {\bf 42}, 141-154 (1978).


\bibitem{Ke69} Kesten H.:  Hitting probabilities of single points for
processes with stationary  independent increments. Memoirs of the
American Mathematical Society, No. 93, American Mathematical
Society, Providence, R.I. (1969).

\bibitem{Luka} Lukacs E.:  Characteristic Functions. 2ed., Griffin, London (1970).



%

\bibitem{PS69} Port S.C., Stone C.J.:  The asymmetric Cauchy
process on the line. Ann. Math. Statist. {\bf 40}, 137-143 (1969).

\bibitem{R77}Rao M.:  On a result of M. Kanda. Z. Wahrsch. verw.
Gebiete {\bf 41}, 35-37 (1977).
%

\bibitem{R88} Rao M.:  Hunt's hypothesis for L\'{e}vy
processes. Proc. Amer. Math. Soc. {\bf 104},  621-624 (1988).

\bibitem{Si77} Silverstein M.L.:  The sector condition implies that
semipolar sets are quasi-polar. Z. Wahrsch. verw. Gebiete {\bf 41},
13-33 (1977).
%


\end{thebibliography}
\end{document}